\documentclass[11pt,twoside,leqno]{aomamlt2e}

\newtheorem{theorem}{Theorem}[section]
\newtheorem{lemma}[theorem]{Lemma}
\theoremstyle{definition}
\newtheorem{definition}[theorem]{Definition}

\newcommand{\abs}[1]{\lvert#1\rvert}

\begin{document}

\title{Discrete differential geometry of triangle tiles and algebra of closed trajectories}

\author{Naoto Morikawa}
\institution{Genocript, 27-22-1015, Sagami-ga-oka 1-chome, Zama-shi, Kanagawa 228-0001 Japan\\
\email{nmorika@f3.dion.ne.jp}}

 \shorttitle{DDG of Triangle Tiles} 

%
%
%


\section{Introduction}
This paper proposes a new mathematical framework that can be applied to 
biological problems such as analysis of the structures of proteins and protein complexes.

In particular, it gives a new method for encoding the three-dimensional structure of a protein into a binary sequence, where proteins are approximated by a folded tetrahedron sequence. The feature of the method is the correspondence between protein folding and ``integration''. And the binary code of a protein is obtained as the ``second derivative" of the native folded structure of the protein. With this method at hand, it becomes possible to describe the structure of a protein without any subjective hierarchical classification.

It also gives a new algebraic framework for describing molecular complexes and their interactions.

First we describe the biological background briefly.

\subsection{The three-dimensional structures of proteins}
Protein is a sequence of amino-acids linked by peptide bonds, where the order of amino-acids are encoded by gene. In nature proteins are folded into a well-defined three-dimensional structure (native state) and this process is called protein folding. Protein folding is reproducible and a protein always folds into the particular naive state rapidly. The native state is characterized by its free energy which is significantly lower than that of alternative structures. (See \cite{B} for more information.)

The functional properties of proteins depend upon their three-dimensional structures. For example, the three-dimensional structure of a protein determines the active site of an enzyme, binding site of a drug, or binding site for another protein. These active and binding sites are key to understanding how proteins interact with other proteins in the cell and how particular molecular targets interact with drugs. 
That is, the knowledge of the three-dimensional structures is crucial to the study of protein function.

Currently, with the growing database of known protein structures, the classification of the structures  plays a central role in understanding the principles of protein structure and function. So far, the structures are classified in a structural hierarchy based on similarity measures. But classification schemes are to some extent subjective and are not amenable to automation. As a result, there are number of different classifications.

\subsection{The structures of protein complexes}
Biological processes, such as signal transmission, cell-fate regulation, transcription, and so on, are not performed by freely diffusing and occasionally colliding proteins. Instead, proteins usually do their jobs by forming structured ensemble of proteins, that is, protein complexes. And one should analyze the structures of protein complexes for a description of biochemical and cellular functions. (See \cite{A} and \cite{S} for more information.)

Frequently, these complexes comprise ten or more subunits. For example, 70S ribosome consists of $52$ proteins and three RNA molecules. But there is no way but rendering a protein complex on a computer graphic system to describe the topological arrangement of its subunits.

\begin{figure}[tb]
\includegraphics{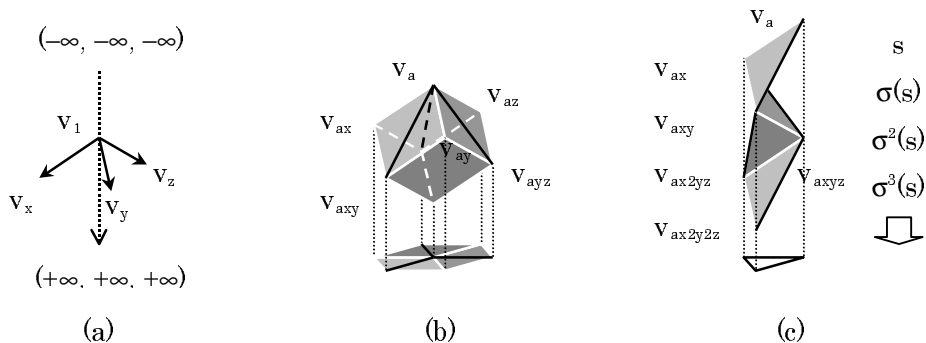}
\caption{Basic ideas. 
(a): Coordinate axes, where $v_{1}=(0,0,0)$, $v_{x}=(1,0,0)$, $v_{y}=(0,1,0)$, and $v_{z}=(0,0,1)$.
(b): Unit cube in $\mathbb{R}^3$ and the projection of the upper faces. 
(c): Slant-tiles over a flat-tile, where $s=v_a[x_{1} x_{2}]$, $\sigma(s)=v_{ax} [x_{2} x_{3}]$, $\sigma^{2}(s)=v_{axy} [x_{3} x_{1}]$, and so on.}
\label{fig1}
\end{figure}

\subsection{Previous works}
Because of the rigidity of covalent bonds between successive amino-acids, the only degree of freedom proteins have are rotations around these covalent bonds. And, traditionally, protein backbone structures are represented by the rotation angles. (See \cite{B} for more information.)

\cite{RS} proposed a representation of the protein backbone structures based on concepts of differential geometry. In their method, a protein is represented as broken lines, where each line corresponds to the virtual-bond between the Ca atoms of successive residues. And they defined the curvature and torsion at each point to describe the native folded structure of a protein.

As for absolute description of protein structure, \cite{T} proposed a set of idealized structures to allow the rigorous and automatic definition of protein topology. And \cite{RF} constructed a measure of similarity to classify protein structures automatically. On the other hand, \cite{H} used the metric matrix distance geometry method to  present a three-dimensional map of the protein structures in which structurally related proteins are represented by spatially adjacent points.

\section{Basic ideas}
The basic idea behind this framework is the following observation: let's consider unit cubes in the $N$-dimensional Euclidean space $\mathbb{R}^N$ and pile them up in the direction from $(+ \infty, + \infty, \cdots, + \infty)$ to $(- \infty, - \infty, \cdots, - \infty)$ and view the resulting surface from $(- \infty, - \infty, \cdots, - \infty)$ (Fig.\ref{fig1}(a)).  If one prints a pattern on the upper faces of each cubes, he/she obtains a drawing made up of the patterns, which defines a flow of $N-1$-dimensional polyhedron tiles of $N$ faces. We call the $N-1$-dimensional polyhedron tiles \textit{triangle tiles} if $N=3$ and \textit{tetrahedron tiles} if $N=4$.

In the case of $N=4$, we obtain a flow of tetrahedrons which we use for encoding of space curves and the structure of proteins. For simplicity, we shall explain the framework in the case of $N=3$, where we use triangle tile sequences for encoding of plane curves and others.

Consider a unit cube in the three-dimensional Euclidean space $\mathbb{R}^3$ specified by the following vertices (Fig.\ref{fig1}(b)): 
\begin{align*}
v_{a}&=(l, m, n) \in \mathbb{Z}^3 \\
v_{ax}&=(l + 1, m, n) \in \mathbb{Z}^3 \\
v_{ay}&=(l, m + 1, n) \in \mathbb{Z}^3 \\
v_{az}&=(l, m, n + 1) \in \mathbb{Z}^3
\end{align*}

And draw lines $\overline{v_{a}v_{axy}}$, $\overline{v_{a}v_{ayz}}$, and $\overline{v_{a}v_{axz}}$ on three faces $v_{a}v_{ax}v_{axy}v_{ay}$, $v_{a}v_{ax}v_{axy}$\\
$v_{ay}$, and $v_{a}v_{ax}v_{axy}v_{a}$ of the cube respectively, where
\begin{align*}
v_{axy}&=(l + 1, m, n) \in \mathbb{Z}^3 \\
v_{ayz}&=(l, m + 1, n + 1) \in \mathbb{Z}^3 \\
v_{axz}&=(l + 1, m, n + 1) \in \mathbb{Z}^3
\end{align*}
Then each face is divided into two slant triangle tiles. For example, triangles $v_{a}v_{ax}v_{axy}$ and $v_{a}v_{ay}v_{axy}$ for $v_{a}v_{ax}v_{axy}v_{ay}$.

\setcounter{figure}{1}

\begin{figure}[tb]
\includegraphics{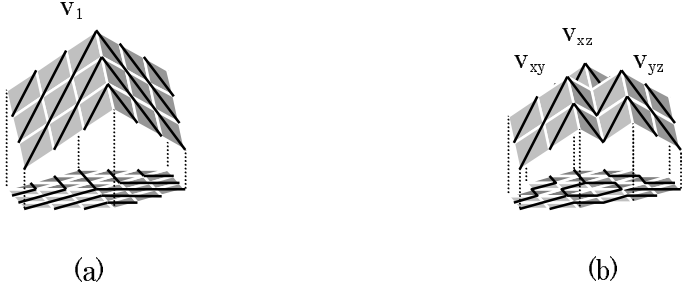}
\caption{Peaks and valleys. 
(a): Drawing defined by a peak $v_1=(0,0,0)$.
(b): Drawing defined by three peaks $v_{xy}=(1,1,0)$, $v_{yz}=(0,1,1)$, and $v_{xz}=(1,0,1)$.}
\label{fig2}
\end{figure}

By piling up these cubes in the direction from $(+ \infty, + \infty, + \infty)$ to $(- \infty,$
$ - \infty, - \infty)$, we obtain ``peaks and valleys'' of cubes with a ``drawing of broken lines" on them. Note that the drawing is uniquely determined by the peaks and divides the surface into sequences of slant triangle tiles. And we obtain a division of $\mathbb{R}^2$ into sequences of flat triangle tiles by ``viewing the surface from $(- \infty, - \infty, - \infty)$''.

For example, the peaks and valleys defined by a peak, say $v_1=(0,0,0)$, form an infinite triangular pyramid, where the top vertex is $v_1$ and the three edges are defined by the positive part of x-, y-, and z-axes (Fig.\ref{fig2}(a)). In this case, each of the three slopes of the cone is divided by a collection of parallel straight lines into infinite sequences of slant triangle tiles.

As another example, let's consider the peaks and valleys defined by three peaks, say $v_{xy}=(1,1,0)$, $v_{yz}=(0,1,1)$, and $v_{xz}=(1,0,1)$ (Fig.\ref{fig2}(b)). Then, the drawing defines a closed sequence of six slant triangle tiles $v_{xy}v_{xyz}v_{x2yz}$, $v_{xz}v_{xyz}v_{x2yz}$, $v_{xz}v_{xyz}v_{xyz2}$, $v_{yz}v_{xyz}v_{xyz2}$, $v_{yz}v_{xyz}v_{xy2z}$, and $v_{xy}v_{xyz}v_{xy2z}$, where $v_{xyz}=(1, 1, 1)$, $v_{x2yz}=(2, 1, 1)$, $v_{xy2z}=(1, 2, 1)$, and $v_{xyz2}=(1, 1, 2)$. (See also Fig.\ref{fig4}(b).)

In the following, we define a vector field on a collection of flat tiles, where a sequence of flat tiles induced by the ``peaks and valleys'' corresponds to a trajectory defined by the vector field.

\section{Cones of three-dimensional lattices}
Here we give the precise definition of the ``peaks and valleys'' of cubes and the drawing on them.

In the first place, we define two kinds of three-dimensional lattices, the standard lattice and its conjugate.
\begin{definition}[$L_3$ and its conjugate $L_3{}^\ast$]
$L_3 := \mathbb{Z}^3$ and $L_3{}^\ast  := \mathbb{Z}^3$, where $L_3$ is embedded in $L_3{}\ast$ and $L_3{}^\ast$ is embedded in $\mathbb{R}^3$ by the following mappings:
\begin{align*}
(l_1, l_2, l_3) \in L_3                   &\mapsto (l_2 + l_3, l_1 + l_3, l_1 + l_2) \in L_3{}^\ast, \\
(m_1, m_2, m_3) \in L_3{}^\ast &\mapsto (m_1, m_2, m_3) \in \mathbb{R}^3.
\end{align*}
\end{definition}

Note that a point of $L_3{}^\ast$ corresponds to a ``fractional point'' of $L_3{}$:
\begin{align*}
(m_1, m_2, m_3) \in L_3{}^\ast \mapsto \bigl( &(m_2 + m_3 - m_1)/2, (m_1 + m_3 - m_2)/2, \\
                                                    & (m_1 + m_2 - m_3)/2 \bigr) \  \in L_3
\end{align*}
if $m_i$s are multiples of two.

Using three indeterminates $x_1$, $x_2$, $x_3$, we obtain the following representation of the lattices (Fig.\ref{fig3}(a)).
\begin{lemma}[Monomial representation of $L_3$ and $L_3{}^\ast$]
\begin{align*}
L_3 &\sim \left\{y_1{}^{l1} y_2{}^{l2} y_3{}^{l3}\ | \   (l_1, l_2, l_3) \in \mathbb{Z}^3  \right\}, \\
L_3{}^\ast &\sim \left\{x_1{}^{m1} x_2{}^{m2} x_3{}^{m3}\ |\  (m_1, m_2, m_3) \in \mathbb{Z}^3 \right\}
\end{align*}
by one-to-one correspondences
\begin{align*}
(l_1, l_2, l_3)  &\sim y_1{}^{l1} y_2{}^{l2} y_3{}^{l3} \\
(m_1, m_2, m_3)  &\sim x_1{}^{m1} x_2{}^{m2} x_3{}^{m3},
\end{align*}
where $y_1 = x_2 x_3$, $y_2 = x_1 x_3$, and $y_3 = x_1 x_2$. 
\end{lemma}

Now let's consider two types of infinite triangular pyramids for each lattice.

\begin{definition}[Standard cones and roofs] For $A \subset L_3{}^\ast$, 
\begin{align*}
Cone\ A &:=\left\{ (l_2 + l_3, l_1 + l_3, l_1 + l_2)  \in \mathbb{R}^3 \ | \  (l_1, l_2, l_3) \in \mathbb{R}^3 \right. \\
& \phantom{:= \{  } \left. \text{and } 
         \exists (a_1, a_2,  a_3)  \in A \text{ s.t. } a_i \le l_i \  (i=1,2,3)\right\}, \\
Roof\ A &:=\left\{ (l_2 + l_3, l_1 + l_3, l_1 + l_2)  \in \mathbb{R}^3 \ | \  (l_1, l_2, l_3) \in \mathbb{R}^3 \right. \\
& \phantom{:= \{  } \left. \text{and }  \exists m \in \mathbb{N} \text{ s.t. } \right.
         (l_2 + l_3, l_1 + m + l_3, l_1 + m + l_2), \\
&\phantom{:= \{ }  (l_2 + m + l_3, l_1 + l_3, l_1 + l_2 + m), (l_2 + l_3 + m, l_1 + l_3 + m, l_1 + l_2) \\
&\phantom{:= \{ } \left. \in Cone\ A \right\}.
\end{align*}
Note that $Cone \  A \subset Roof \  A$ for any $A  \subset L_3{}^\ast$. $Cone\ A$ is called \textit{standard cone} and $Roof\ A$ is called \textit{standard roof}.
\end{definition}

\begin{definition}[Conjugate cones and roofs] For $A \subset L_3{}^\ast$, we define
\begin{align*}         
&Cone^\ast  A:=\left\{ ( l_1, l_2,  l_3) \in \mathbb{R}^3  \ | \ 
         \exists (a_1, a_2,  a_3)  \in A \text{ s.t. } a_i \le l_i \  (i=1,2,3)\right\}, \\
&Roof^\ast  A:=\left\{(l_1, l_2,  l_3) \in \mathbb{R}^3  \ | \  \exists m \in \mathbb{N} \text{ s.t. } 
         (l_1 + m, l_2,  l_3), (l_1, l_2+m,  l_3), \right. \\
&\phantom{Roof\ A:=\{(l_1, l_2,  l_3) \in L_3 \ | \  }  \left.       (l_1, l_2,  l_3 + m) \in Cone^\ast  A \right\}.
\end{align*}
Note that $Cone^\ast A \subset Roof^\ast A$ for any $A  \subset L_3{}^\ast$. $Cone^\ast A$ is called \textit{conjugate cone} and $Roof^\ast A$ is called \textit{conjugate roof}.
\end{definition}

\setcounter{figure}{2}
\begin{figure}[tb]
\includegraphics{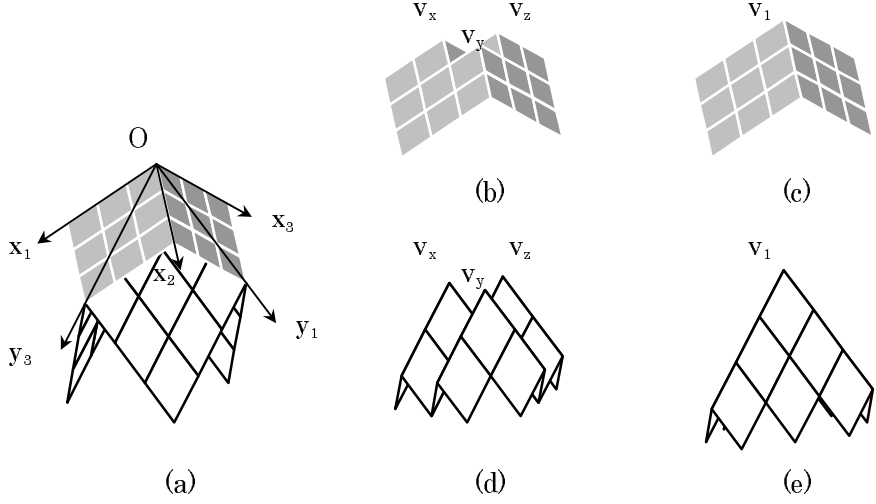}
\caption{Lattices and cones. 
(a): Two lattices $L_3$ (white) and $L_3{}^\ast$ (gray).
(b): $Cone^\ast  \{ v_x, v_y, v_z \}$.
(c): $Roof^\ast  \{ v_x, v_y, v_z \} = Cone^\ast  \{ v_1 \}$.
(d): $Roof\  \{ v_x, v_y, v_z \} = Cone\  \{ v_x, v_y, v_z \}$.
(e): $Cone\  \{ v_1 \}$.}
\label{fig3}
\end{figure}

For example, $Roof \   \{ v_x, v_y, v_z \} = Cone \  \{ v_x, v_y, v_z \}$ and $Roof^\ast  \{ v_x, v_y, v_z \}$
$ = Cone^\ast  \{ v_1 \}$, where $v_{x}=(1, 0, 0)$, $v_{y}=(0, 1, 0)$, and $v_{z}=(0, 0, 1)$  (Fig.\ref{fig3}).

Using indeterminates $x_1$, $x_2$, $x_3$, we obtain the following.
\begin{lemma}[Monomial representation] For $A \subset L_3{}^\ast$, 
\begin{align*}
Cone\ A &\sim \left\{y_1{}^{l1} y_2{}^{l2} y_3{}^{l3} \ | \  ( l_1, l_2,  l_3) \in \mathbb{R}^3 \text{ and } \right. \\
& \phantom{\sim \{y_1{}^{l1} y_2{}^{l2} y_3{}^{l3} \ | \  } \left.
         \exists (a_1, a_2,  a_3)  \in A \text{ s.t. } a_i \le l_i \  (i=1,2,3) \right\}, \\
Roof\ A &\sim \left\{ y_1{}^{l1} y_2{}^{l2} y_3{}^{l3} \ | \  ( l_1, l_2,  l_3) \in \mathbb{R}^3 \text{ and }
\exists m \in \mathbb{N} \text{ s.t. }  \right. \\
& \phantom{\sim \{y_1{}^{l1} y_2{}^{l2} y_3{}^{l3} \ | \  } \left.
         y_1{}^{l1} y_2{}^{l2} y_3{}^{l3} y_i{}^m  \in Cone\ A \  (i=1,2,3)  \right\}, \\
Cone^\ast  A &\sim \left\{x_1{}^{l1} x_2{}^{l2} x_3{}^{l3} \ | \  ( l_1, l_2,  l_3) \in \mathbb{R}^3 \text{ and } \right. 
\end{align*}
\begin{align*}
& \phantom{\sim \{y_1{}^{l1} y_2{}^{l2} y_3{}^{l3} \ | \  } \left.
         \exists (a_1, a_2,  a_3)  \in A \text{ s.t. } a_i \le l_i \  (i=1,2,3) \right\}, \\
Roof^\ast  A &\sim \left\{ x_1{}^{l1} x_2{}^{l2} x_3{}^{l3} \ | \  ( l_1, l_2,  l_3) \in \mathbb{R}^3 \text{ and } \exists m \in \mathbb{N} \text{ s.t. }  \right. \\
& \phantom{\sim \{y_1{}^{l1} y_2{}^{l2} y_3{}^{l3} \ | \  } \left.
         x_1{}^{l1} x_2{}^{l2} x_3{}^{l3} x_i{}^m  \in Cone^\ast  A \  (i=1,2,3)  \right\}.
\end{align*}
\end{lemma}

We denote the ``peaks'' of a cone $w$ by $g(w)$. That is, $g(w) \subset L_3$ (or $g(w) \subset L_3{}^\ast$) is the minimal system of elements of $w$ which satisfies $w=Cone\ g(w)$ (or  $w=Cone^\ast g(w)$). In particular, $Roof\ A = Cone\ g(Roof\ A)$ for $A \subset L_3$ and $Roof^\ast A = Cone^\ast g(Roof ^\ast A)$ for $A \subset L_3{}^\ast$. 
The boundary surface $dw \subset \mathbb{R}^3$ of a conjugate cone (or roof) $w$ is given by:
\[
dw :=  \left\{ (l_1, l_2,  l_3) \in w \  | \  \max_{(a_1, a_2,  a_3) \in g(w)} \left\{  \min \{l_1 - a_1, l_2 - a_2, l_3 - a_3\} \right\} = 0 \right\}.
\]
The boundary surface $dw \subset \mathbb{R}^3$ of a standard cone (or roof) $w$ is given by:
\begin{align*}
dw := & \left\{ (l_2 + l_3, l_1 + l_3, l_1 + l_2)  \in w \  | \   ( l_1, l_2,  l_3) \in \mathbb{R}^3 \text{ and } \right. \\
&\phantom{\{ \quad  } \left. \max_{(a_1, a_2,  a_3) \in g(w)} \left\{  \min \{l_1 - a_1, l_2 - a_2, l_3 - a_3\} \right\} = 0 \right\}.
\end{align*}

Then, the ``peaks and valleys'' of cubes are defined as the boundary surface of a conjugate cone. For example, one defined by a peak $v \in \mathbb{Z}^3$ is $d(Cone^\ast  \{v\})$ and another one defined by three peaks $v_{xy}$, $v_{yz}$, and $v_{xz} \in \mathbb{Z}^3$ is  $d(Cone^\ast \{ v_{xy}, v_{yz}, v_{xz} \})$.

\section{Triangle tiles}
Next we give the definition of the triangle tiles.
For $v_1$, $v_2$, and $v_3 \in \mathbb{R}^3$, we denote the convex hull of them by $conv[v_1, v_2, v_3]$, that is,
\begin{align*}
conv\left[v_1, v_2, v_3 \right]:= & \left\{ \sum_{i =1,2,3}   \lambda_i v_i  \in \mathbb{R}^3 \ |\   \lambda_i \in \mathbb{R} \ (i=1,2,3) \text{ s.t. }   \right. \\
& \phantom{\{ \sum_{i = 1,2,3}  \lambda_i v_i \in \mathbb{R}^3 \ |\  \  } \left. \sum_{i =1,2, 3} \lambda_i=1  \text{ and } \lambda_i \ge 0 \  (i=1,2,3) \right\}.
\end{align*}

\begin{definition}[Slant triangle tiles] For $v_a = (l, m, n) \subset L_3{}^\ast$, we define six \textit{slant triangle tiles} (Fig.\ref{fig1}(b)):
\begin{align*}
v_a\left[x_1 x_2 \right] & := conv\left[v_a, v_{ax}, v_{axy} \right], \\
v_a\left[x_1 x_3 \right] & := conv\left[v_a, v_{ax}, v_{axz} \right], \\
v_a\left[x_2 x_1 \right] & := conv\left[v_a, v_{ay}, v_{ayx} \right], \\
v_a\left[x_2 x_3 \right] & := conv\left[v_a, v_{ay}, v_{ayz} \right], \\
v_a\left[x_3 x_1 \right] & := conv\left[v_a, v_{az}, v_{azx} \right], \\
v_a\left[x_3 x_2 \right] & := conv\left[v_a, v_{az}, v_{azy} \right], 
\end{align*}
where  $v_{ax}=(l+1, m, n)$, $v_{ay}=(l, m+1, n)$, and $v_{az}=(l, m, n+1)$, $v_{axy}=(l+1, m+1, n)$, $v_{ayz}=(l, m+1, n+1)$, $v_{axy}=(l+1, m+1, n)$, and $v_{axz}=(l+1, m, n+1)  \in L_3{}^\ast$.
We denote  the collection of all slant triangle tiles by $S_3$:
\[
S_3 := \left\{ a\left[x_{\rho(1)}  x_{\rho(2)}\right] \ |\   a \in L_3{}^\ast, \  \rho \in \mathit{Sym_3} \right\},
\]
where $\mathit{Sym_3}$ is the collection of all permutations of three letters, 
 that is, the symmetric group on three letters.
\end{definition}

Finally, we give the definition of the ``view from $(- \infty, - \infty, - \infty)$''. To identify slant triangle tiles in the direction from $(+ \infty, + \infty, + \infty)$ to $(- \infty, - \infty,$
$ - \infty)$, we consider  ``shift operator'' $\sigma$ on $S_3$ (Fig.\ref{fig1}(c)):
\[
\sigma\left(a\left[x_{\rho(1)} \cdots x_{\rho(N-1)}\right]\right) 
:= ax_{\rho(1)}\left[x_{\rho(2)} \cdots x_{\rho(N)}\right], 
\]
where $ax_1$, $ax_2$, and $ax_3$ denote the point of $(l+1,m,n)$, $ (l,m+1,n)$, and $(l,m,n+1) \in \mathbb{Z}^3$ respectively for $a = (l,m,n)$. 
Note that shift operator $\sigma$ induces an equivalence relation on $S_3$:
\[
s_1 \sim_\sigma s_2  \textit{   if and only if   } \exists k \in \mathbb{Z}  \text{ s.t. } \sigma^k(s_1) = s_2,
\]
for $s_1$, $s_2 \in S_3$.

\begin{definition}[Flat triangle tiles] A \textit{flat triangle tile} is an $\sim_\sigma$-
 equivalence class of a slant triangle tile and denoted by $s \mod \sigma$ for some $s \in S_3$. In particular, the collection $B_3$ of all \textit{flat triangle tiles} is the quotient of $S_3$ by $\sim_\sigma$: $B_3:= S_3 / \sim_\sigma$.
\end{definition}

For example, let's consider the cone $w = Cone^\ast \{ v_{xy}, v_{yz}, v_{xz} \}$ (Fig.\ref{fig2}(b)). Then, the area of $dw$ enclosed by $dw \cap d(Cone \ \{v_1\})$ consists of six slant triangle tiles, $v_{xy}[x_3x_1]$, $v_{xz}[x_2 x_1]$, $v_{xz}[x_2 x_3]$, $v_{yz}[x_1 x_3]$, $v_{yz}[x_1 x_2]$, and $v_{xy}[x_3 x_2] \in S_3$ (Fig.\ref{fig4}(b)). And the corresponding six flat triangle tiles form a regular hexagon, that is the ``view of the area from $(- \infty, - \infty, - \infty)$''.

\section{Surface decomposition}
Now we give the definition of  the ``drawing'' on the surface of a conjugate cone, which divides the surface into sequences of slant triangle tiles.

Firstly, we define a decomposition of slant tiles contained in the surface of a conjugate cone.
\begin{definition}[Decomposition by a standard cone] For  a conjugate\\
 cone $w_1$ and a standard cone $w_2$, set 
\begin{align*}
In(w_1, w_2) &:= \left\{ s \in S_3 \ | \  s \subset dw_1 \cap w_2   \right\}, \\
Out(w_1, w_2) &:= \left\{ s \in S_3 \ | \  s \subset dw_1 \setminus w_2   \right\}, \\
Bd(w_1, w_2) &:= \left\{ s \in dw_1 \ | \  s \notin In(w_1, w_2) \cup Out(w_1, w_2)   \right\}.
\end{align*}
$w_2$ is called \textit{consistent} with $w_1$ if $Bd(w_1. w_2) = \emptyset$.
\end{definition}

\setcounter{figure}{3}
\begin{figure}[tb]
\includegraphics{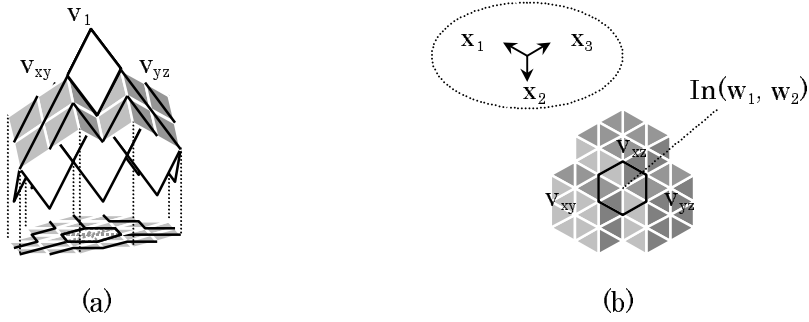}
\caption{Surface decomposition.
(a): $w_1 = Cone^\ast \{ v_{xy}, v_{yz}, v_{xz} \}$ (gray) and $w_2 = Cone \  \{v_1 \}$ (white).
(b): $In(w_1, w_2)$.}
\label{fig4}
\end{figure}

Then, for a conjugate cone $w$, the ``drawing'' on $dw$ is defined as the collection of all intersections between $dw$ and the surfaces of the standard cones which are consistent with $w$:
\[
\{ dw \cap d(Cone\  A) \  | \  A \subset L_3  \text{ and } Bd(w, Cone\ A) = \emptyset \}  \subset dw.
\]

For example, we obtain
\begin{align*}
In(w_1, w_2) &= \left\{ v_{xy}[x_3x_1], v_{xz}[x_2 x_1], v_{xz}[x_2 x_3], v_{yz}[x_1 x_3], v_{yz}[x_1 x_2], v_{xy}[x_3 x_2]    \right\}, \\
Bd(w_1, w_2) &= \emptyset
\end{align*}
for $w_1 = Cone^\ast \{ v_{xy}, v_{yz}, v_{xz} \}$ and $w_2 = Cone \  \{v_1 \}$ (Fig.\ref{fig4}(a)). Then, intersection $dw_1 \cap dw_2$ gives the boundary of deformed hexagon  $v_{xy}$$v_{x2yz}$
$v_{xz}$$v_{xyz2}$$v_{yz}$
$v_{xy2z}$ (Fig.\ref{fig4}(b)).

\section{Differential geometry of triangle tiles}

\subsection{Differential structure on $B_3$}
To define ``tangent bundle'' on $B_3$, we consider the equivalence relation on $S_3$ induced by  $\sigma^3$:
\[
s_1 \sim_{\sigma 3} s_2  \textit{   if and only if   } \exists k \in \mathbb{Z}  \text{ s.t. } \sigma^{3k}(s_1) = s_2,
\]
for $s_1$, $s_2 \in S_3$. We denote an $\sim_{\sigma 3}$-equivalence class by $s \mod \sigma^3$ for some $s \in S_3$.

\begin{definition}[Tangent bundle on $B_3$]
\textit{Tangent bundle} $T[B_3]$ on $B_3$ is the quotient of $S_3$ by $\sim_{\sigma 3}$:
\begin{align*}
&T[B_3]:= S_3/ \sim_{\sigma 3}, \\ 
&\pi: T[B_3] \to B_3, \  \pi\left(s \mod \sigma^3 \right):= s \mod \sigma.
\end{align*}
$T[B_3]$ is identified with $B_3 \times \{x_1 x_2, x_2 x_3, x_1 x_3 \}$ by one-to-one correspondence
\[
s \mod \sigma^3  \sim ( s \mod \sigma, Ds),
\]
where the \textit{gradient} $Ds$ of $s \in S_3$ is defined by
\[
D a\left[x_{\rho(1)} x_{\rho(2)}\right] :=x_{\rho(1)} x_{\rho(2)}.
\]
\end{definition}

For each element of $T[B_3]$, we assign a local trajectory. And we shall obtain a flow on $B_3$ by patching these local trajectories together. 

We start with the definition of the local trajectories on $S_3$.
\begin{definition}[Local trajectories on $S_3$]
For $a\left[x_{\rho(1)} x_{\rho(2)}\right]$ $\in S_3$, we assign the following four trajectories of length three  (Fig.\ref{fig5}(a)):
\begin{align*}
& \left\{ a/x_{\rho(2)}  \left[x_{\rho(2)} x_{\rho(1)}\right],\    a\left[x_{\rho(1)} x_{\rho(2)}\right] ,\quad ax_{\rho(1)}/x_{\rho(3)} \left[x_{\rho(3)} x_{\rho(2)}\right] \right\} \\
& \left\{ a\left[x_{\rho(1)} x_{\rho(3)}\right],\quad \quad \quad  a\left[x_{\rho(1)} x_{\rho(2)}\right] ,\quad  ax_{\rho(1)}/x_{\rho(3)} \left[x_{\rho(3)} x_{\rho(2)}\right] \right\} \\
& \left\{ a/x_{\rho(2)}  \left[x_{\rho(2)} x_{\rho(1)}\right],\    a\left[x_{\rho(1)} x_{\rho(2)}\right] ,\quad  ax_{\rho(1)}\left[x_{\rho(2)} x_{\rho(1)}\right] \right\} \\
& \left\{ a\left[x_{\rho(1)} x_{\rho(3)}\right],\quad \quad \quad  a\left[x_{\rho(1)} x_{\rho(2)}\right] ,\quad  ax_{\rho(1)}\left[x_{\rho(2)} x_{\rho(1)}\right] \right\}
\end{align*}
A \textit{flow} on $S_3$ is obtained by patching one of these local trajectories together. 
\end{definition}

\setcounter{figure}{4}
\begin{figure}[tb]
\includegraphics{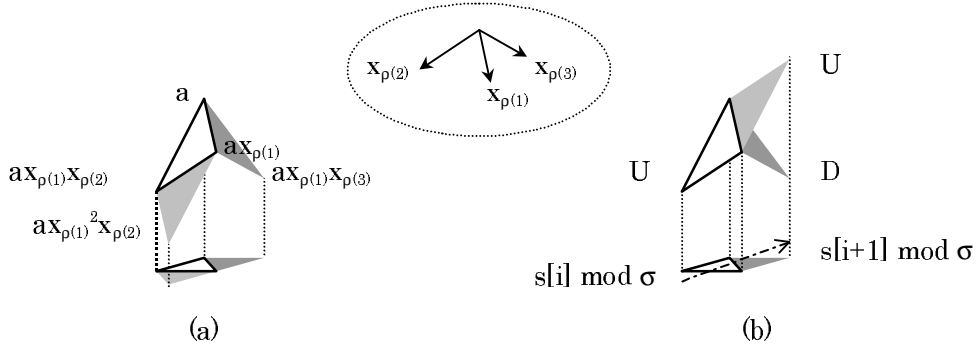}
\caption{Differential structure.
(a): Local trajectories specified by $a[x_{\rho(1)} x_{\rho(2)}]$ $\in S_3$.
(b): The second derivative along $\{s[i]\}$.}
\label{fig5}
\end{figure}

Note that they are mapped  onto the same local trajectory on $B_3$ by $\pi$.

\begin{lemma}[Local trajectories on $B_3$]
The local trajectories on $S_3$ at $a\left[x_{\rho(1)} x_{\rho(2)}\right]$ $\in S_3$ induce the following local trajectory on $B_3$ at $a\left[x_{\rho(1)} x_{\rho(2)}\right]  \mod \sigma$:
\[
\left\{ a\left[x_{\rho(1)} x_{\rho(3)}\right]  \mod \sigma,
\  a\left[x_{\rho(1)} x_{\rho(2)}\right]  \mod \sigma,
\  ax_{\rho(1)}\left[x_{\rho(2)} x_{\rho(1)}\right]  \mod \sigma \right\}
\]
\end{lemma}

The surface of a conjugate cone defines a flow on $S_3$ and we have the following.

\begin{lemma}Let $\{s[i] \  | \  1 \le i \le n\} \subset S_3$ be a closed trajectory defined by a conjugate cone $w$. Then, there exist standard roofs $w_1$ and $w_2$ s.t. they are consistent with $w$ and
\[
\{s[i]  \  | \  1 \le i \le n \}  = In(w, w_1) \setminus In(w, w_2).
\]
\end{lemma}

\Proof
Suppose that $s[i] = a_i [x_{i1}x_{i2}]$ ($1 \le i \le n$) and set 
\begin{align*}
w_1 &= Roof \  \{ a_i \ | \  1 \le i \le n \}, \\
w_2 &= Roof \ \{ a \ |\   \exists s \in B  \text{  s.t. }  s  = a[x_{\rho(1)} x_{\rho(2)}]  \text{ for some } \rho \in Sym_3 \},
\end{align*}
where $B = \{ s \in In(w, w_1) \ | \  s  \neq s[i] \ ( 1 \le i \le n)\}$.
\Endproof

Since there are only two choices of succeeding slant tiles for each direction at $a\left[x_{\rho(1)} x_{\rho(2)}\right]$ $\in S_3$:
\[
a/x_{\rho(2)}  \left[x_{\rho(2)} x_{\rho(1)}\right]  \quad \text{ or }\quad a\left[x_{\rho(1)} x_{\rho(3)}\right]
\]
for one direction and
\[
ax_{\rho(1)}/x_{\rho(3)} \left[x_{\rho(3)} x_{\rho(2)}\right] \quad  \text{ or } \quad  ax_{\rho(1)}\left[x_{\rho(2)} x_{\rho(1)}\right]
\]
for the other, we obtain the following definition of the ``second derivative'' along a trajectory on $S_3$.

\begin{definition}[The second derivative] Let $\{ s[i] \}$ be a trajectory on $S_3$. 
The \textit{second derivative} $D^2s[i]$ along $\{s[i]\}$ is a $\{ U, D \}$-valued function defined by
\[
D^2(s[i+1]):=  \begin{cases}
                   D^2(s[i])  \quad  \text{if  $D(s[i+1])=D(s[i])$}, \\
                   - D^2(s[i])  \quad  \text{else},
                   \end{cases}
\]
where $- D:= U$ and $-U:=D$ (Fig.\ref{fig5}(b)). 
In other words, the value of the second derivative is negated when the gradient of a trajectory changes.
\end{definition}

Note that we can decode the shape of a trajectory on $S_3$ by the second derivative along it, i.e., by an $U/D$ sequence. See below for an example, an encoding of the shape of a regular hexagon.

\subsection{Vector field on $B_3$ induced by a cone}
For  a conjugate cone, $dw$ specifies a unique slant triangle tile over 
each point of $B_3$. In other words, $dw$ specifies a section of $S_3$ over $B_3$, which we denote by $\Gamma_w$:
\[
\Gamma_w(t):= \text{ the unique slant tile $s \subset dw$ \quad   s.t.  \quad  $t=s \mod \sigma$}
\]
for $t \in B_3$. Then, $\Gamma_w$ induces a vector field over $B_3$.

\begin{definition}[The vector field induced by a conjugate cone $w$]
\[
 X_w(s \mod \sigma ):=D\Gamma_w(s \mod \sigma ) \in \{ x_1 x_2, x_2 x_3, x_1 x_3 \}. 
\]
\end{definition}

And any vector field on $B_3$ is locally defined by a conjugate cone. In general, we need more than one cone to specify a vector field because of overlaps between slant triangle tiles as in the case of local charts of a manifold of conventional differential geometry.

Let $\{ t[i] \} \subset B_3$ be a trajectory defined by a vector field $X_w$. Then the second derivative $D^2 \Gamma_w$ of $\Gamma_w$ along $\{t[i]\}$ is given by
\[
D^2 \Gamma_w(t[i+1]) =  \begin{cases}
                   D^2 \Gamma_w(t[i])  \quad  \text{if  $X_w(t[i+1])=X_w(t[i])$}, \\
                   - D^2 \Gamma_w(t[i])  \quad  \text{else}.
                   \end{cases}
\]

For example, let's consider the conjugate cone $w = Cone^\ast \{ v_{xy}, v_{yz}, v_{xz} \}$ again (Fig.\ref{fig4}). We have seen that the area of $dw$ enclosed by $dw \cap d(Cone \ \{v_1\})$ consists of six slant triangle tiles and mapped onto a regular hexagon by $\pi$. Actually they form a closed trajectory $\{t[i] \ | \  1 \le i \le 6\}$ of $B_3$, where
\begin{align*}
&\Gamma_w(t[1] ) = v_{xy}[x_3x_1],  \quad X_w(t[1] )  = x_1 x_3, \quad D^2 \Gamma_w(t[1]) = D, \\
&\Gamma_w(t[2] ) = v_{xz}[x_2 x_1],  \quad X_w(t[2] )  = x_1 x_2, \quad D^2 \Gamma_w(t[2]) = U, \\
&\Gamma_w(t[3] ) = v_{xz}[x_2 x_3],  \quad X_w(t[3] )  = x_2 x_3, \quad D^2 \Gamma_w(t[3])  = D, \\
&\Gamma_w(t[4] ) = v_{yz}[x_1 x_3],  \quad X_w(t[4] )  = x_1 x_3, \quad D^2 \Gamma_w(t[4])  = U, \\
&\Gamma_w(t[5] ) = v_{yz}[x_1 x_2],  \quad X_w(t[5] )  = x_1 x_2, \quad D^2 \Gamma_w(t[5])  = D, \\
&\Gamma_w(t[6] ) = v_{xy}[x_3 x_2],  \quad X_w(t[6] )  = x_2 x_3, \quad D^2 \Gamma_w(t[6])  = U.
\end{align*}

As a result, we obtain a binary code of the shape of a regular hexagon:
\[
D-U-D-U-D-U.
\]
Note that, if we had set $D^2 \Gamma_w(t[1])$ to $U$, the binary code would be negated.

On the other hand, we can compute a trajectory which is encoded by a given sequence of $U$ and $D$. For example, let's consider the binary sequence of 
\[
D-U-D-U-D-D
\]
and construct a section $\Gamma$ of $S_3$ over $B_3$ which induces a trajectory whose second derivative is given by the sequence. Suppose $\Gamma(t[1])= v_{xy}[x_3x_1]$ and $t[2]=v_{xz}[x_2 x_1] \mod \sigma$. Then, we obtain
\begin{align*}
&\Gamma(t[1] ) = v_{xy}[x_3x_1],  \quad \ D\Gamma(t[1] )  = x_1 x_3, \\
&\Gamma(t[2] ) = v_{xz}[x_2 x_1],  \quad \ D\Gamma(t[2] )  = x_1 x_2, \\
&\Gamma(t[3] ) = v_{xz}[x_2 x_3], \quad \ D\Gamma(t[3] )  = x_2 x_3, \\
&\Gamma(t[4] ) = v_{yz}[x_1 x_3],  \quad \ D\Gamma(t[4] )  = x_1 x_3, \\
&\Gamma(t[5] ) = v_{yz}[x_1 x_2],   \quad \ D\Gamma(t[5] )  = x_1 x_2, \\
&\Gamma(t[6] ) = v_{xyz}[x_2 x_1], \quad D\Gamma(t[6] )  = x_1 x_2.
\end{align*}
The first five tiles $t[1], t[2], \cdots, t[5]$ correspond to the vector field induced by 
$Cone^\ast  \left\{  v_{xy}, v_{yz}, v_{xz}  \right\}$  and the last five tiles $t[2], t[3], \cdots, t[6]$  correspond to the vector field induced by $Cone^\ast  \left\{  v_{yz}, v_{xz}  \right\}$. That is, we need two local vector fields to cover the trajectory.

\section{Algebra of closed trajectories}

\subsection{Addition of roofs}
Finally we define addition of roofs to consider an analog of protein complexes.

\begin{definition}[Addition of conjugate roofs]
For $A, B \subset L_3{}^\ast$, 
\[
Roof^\ast A + Roof^\ast B := Roof^\ast A\cup B.
\]
\end{definition}

For example, 
\begin{align*}
Roof^\ast \{ v_{x} \} + Roof^\ast \{ v_{y} \} + Roof^\ast \{ v_{z} \} &= Roof^\ast \{ v_1 \}, \\
Roof^\ast \{ v_{xy} \} + Roof^\ast \{ v_{yz} \} + Roof^\ast \{ v_{xz} \} &= Roof^\ast \{ v_{xy}, v_{yz}, v_{xz} \}.
\end{align*}

And we set $\abs{w} := \pi(In(w, Roof \  g(w)) \subset B_3$ for a roof $w$ of $L_3{}^\ast$, which we call the \textit{norm} of $w$.  For example,
\begin{align*}
\left|Roof^\ast \{ v_{1} \} \right| &= \left|Roof^\ast \{ v_{x} \} \right| = \cdots = \left|Roof^\ast \{ v_{xz} \} \right| =\emptyset, \\
\left|Roof^\ast \{ v_{xy}, v_{yz}, v_{xz} \} \right| &= 
\left\{ v_{xy}[x_3x_1] \mod \sigma, \quad v_{xz}[x_2 x_1] \mod \sigma, \right.\\
&\phantom{= \quad } \  v_{xz}[x_2 x_3] \mod \sigma, \quad v_{yz}[x_1 x_3] \mod \sigma, \\
&\phantom{= \quad }  \left. v_{yz}[x_1 x_2] \mod \sigma, \quad v_{xy}[x_3 x_2] \mod \sigma    \right\}.
\end{align*}
(Note that $Roof^\ast \{ v_{xy}, v_{yz}, v_{xz} \} = Cone^\ast \{ v_{xy}, v_{yz}, v_{xz} \}$ and see Fig.\ref{fig4}(b).)

In particular, we obtain the following.

\begin{lemma} For a finite set $\{ w_i \}$ of conjugate roofs, $\bigcup_i \left| w_i \right| \subset \left| \sum_i w_i \right|$.
\end{lemma}

For another example, set
\begin{align*}
w_1 &= Roof^\ast \{ v_{xy}, v_{yz}, v_{xz} \}, \\
w_2 &= Roof^\ast \{ v_{xy2z-1}, v_{y2}, v_{xy} \}, \\
w_3 &= Roof^\ast \{ v_{y2}, v_{x-1y2z}, v_{yz} \},
\end{align*}
where $v_{xy2z-1}=(1,2,-1)$, $v_{y2}=(0,2,0)$, and $v_{x-1y2z}=(-1,2,1)$ (Fig.\ref{fig6}(a)). Note that $\abs{w_i}$ ($i=1,2,3$) is a closed trajectory of length six and forms a regular hexagon. Then, we obtain
\begin{align*}
w_1 + w_2 + w_3 &= Roof^\ast \{ v_{xy}, v_{xy-1z2}, v_{x-1yz2}, v_z \}, \\
\abs{w_1} \cup \abs{w_2} \cup \abs{w_3} &= \abs{w_1+w_2+w_3}.
\end{align*}
That is, three trajectories $\abs{w_i}$ ($i=1,2,3$) of length six are fused into a trajectory of length $18$ (Fig.\ref{fig6}(b)). 
In other words, $\abs{w_1+w_2+w_3}$ is an analog of protein complexes made up of three subunits.

\setcounter{figure}{5}
\begin{figure}[tb]
\includegraphics{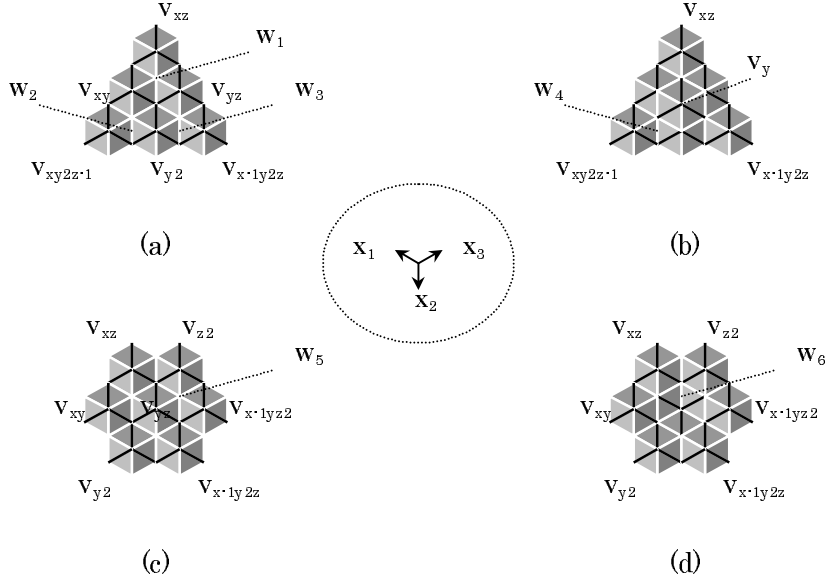}
\caption{Algebra of closed trajectories.
(a): The closed trajectories defined by $w_1$, $w_2$, and $w_3$.
(b): The closed trajectory defined by $w_4$.
(c): The closed trajectories defined by $w_2 + w_3 + w_5$.
(d): The closed trajectory defined by cone $w_6$.}
\label{fig6}
\end{figure}

In general, a conjugate roof is not associated with a single trajectory. But there may exist a conjugate cone which defines a single closed trajectory which sweeps all the closed trajectories of the roof. For example, consider $w_1 + w_3 + w_5$, where $w_5 = Roof^\ast \{ v_{yz}, v_{x-1yz2}, v_{z2} \}$ (Fig.\ref{fig6}(c)). It  is associated with three closed trajectories and
\[
w_1 + w_3 + w_5 = Roof^\ast \{ v_{xz}, v_{xy}, v_{y2}, v_{x-1y2z}, v_{x-1yz2}, v_{z2} \}
\]
Then, $w_6 = Cone^\ast \{ v_{xz}, v_{xy}, v_{y2}, v_{x-1y2z}, v_{x-1yz2}, v_{z2} \}$ defines a single closed trajectory of $\abs{w_1+w_3+w_5}$ (Fig.\ref{fig6}(d)).

\section{Conclusion}
The basic idea of a new mathematical framework that can be applied to biological problems such as analysis of the structure of proteins and protein complexes has been described.

 If we consider the case of $N=4$, we would obtain a description of proteins as trajectories of tetrahedron tiles, where the structure of a protein is encoded into a sequence of $U$ and $D$  (\cite{M}) . And we could extract static structual information of proteins directly from genes by comparing $U/D$ sequences with the corresponding genes that are sequences of four letters A, T, G, and C.

 Moreover, by assigning roofs of $L_3{}^\ast$ to proteins and their complexes, we could describe formations of protein complexes algebraically.

\end{document}